\magnification 1250
\pretolerance=500 \tolerance=1000  \brokenpenalty=5000
\mathcode`A="7041 \mathcode`B="7042 \mathcode`C="7043
\mathcode`D="7044 \mathcode`E="7045 \mathcode`F="7046
\mathcode`G="7047 \mathcode`H="7048 \mathcode`I="7049
\mathcode`J="704A \mathcode`K="704B \mathcode`L="704C
\mathcode`M="704D \mathcode`N="704E \mathcode`O="704F
\mathcode`P="7050 \mathcode`Q="7051 \mathcode`R="7052
\mathcode`S="7053 \mathcode`T="7054 \mathcode`U="7055
\mathcode`V="7056 \mathcode`W="7057 \mathcode`X="7058
\mathcode`Y="7059 \mathcode`Z="705A
\def\spacedmath#1{\def\packedmath##1${\bgroup\mathsurround 
=0pt##1\egroup$}\mathsurround#1
\everymath={\packedmath}\everydisplay={\mathsurround=0pt}}
 \spacedmath{2pt}
\def\mono{\lhook\joinrel\mathrel{\longrightarrow}}
\def\iso{\vbox{\hbox to .8cm{\hfill{$\scriptstyle\sim$}\hfill}
\nointerlineskip\hbox to .8cm{{\hfill$\longrightarrow $\hfill}} }}
\def\sdir_#1^#2{\mathrel{\mathop{\kern0pt\oplus}\limits_{#1}^{#2}}}
\def\pprod_#1^#2{\raise
2pt \hbox{$\mathrel{\scriptstyle\mathop{\kern0pt\prod}\limits_{#1}^{#2}}$}}

\font\san=cmssdc10\font\gragrec=cmmib10
\def\mug{\hbox{\gragrec \char22}}

\def\sym{\hbox{\san \char83}}

\def\pc#1{\tenrm#1\sevenrm}
\def\tx{\kern-1.5pt -}
\def\cqfd{\kern 2truemm\unskip\penalty 500\vrule height 4pt depth 0pt width
4pt\medbreak} 
\def\virg{\raise
.4ex\hbox{,}}
\def\ind{\par\hskip 1truecm\relax}

\def\rond{\kern 1pt{\scriptstyle\circ}\kern 1pt}
\def\iso{\mathrel{\mathop{\kern 0pt\longrightarrow }\limits^{\sim}}}

\def\Hom{\mathop{\rm Hom}\nolimits}
\def\Aut{\mathop{\rm Aut}\nolimits}

\frenchspacing

\font\eightbboard=msbm8\font\tenbboard=msbm10
\font\sevenbboard=msbm7
\newfam\bboardfam
\textfont\bboardfam=\tenbboard
\scriptfont\bboardfam=\sevenbboard
 \def\bb{\fam\bboardfam}

\font\igr=cmmib7
\def\nn{\hbox{\igr\char110}}
\font\eightrm=cmr8         \font\eighti=cmmi8
\font\eightsy=cmsy8        \font\eightbf=cmbx8
\font\eighttt=cmtt8        \font\eightit=cmti8
\font\eightsl=cmsl8        \font\sixrm=cmr6
\font\sixi=cmmi6           \font\sixsy=cmsy6 

\font\sixbf=cmbx6\catcode`\@=11
\def\eightpoint{%
  \textfont0=\eightrm \scriptfont0=\sixrm \scriptscriptfont0=\fiverm
  \def\rm{\fam\z@\eightrm}%
  \textfont1=\eighti  \scriptfont1=\sixi  \scriptscriptfont1=\fivei
  \textfont2=\eightsy \scriptfont2=\sixsy \scriptscriptfont2=\fivesy
  \textfont\itfam=\eightit
  \def\it{\fam\itfam\eightit}%
 \textfont\bboardfam=\eightbboard
  \def\bb{\fam\bboardfam}
  \textfont\slfam=\eightsl
  \def\sl{\fam\slfam\eightsl}%
  \textfont\bffam=\eightbf \scriptfont\bffam=\sixbf
  \scriptscriptfont\bffam=\fivebf
  \def\bf{\fam\bffam\eightbf}%
  \textfont\ttfam=\eighttt
  \def\tt{\fam\ttfam\eighttt}%
  \abovedisplayskip=9pt plus 3pt minus 9pt
  \belowdisplayskip=\abovedisplayskip
  \abovedisplayshortskip=0pt plus 3pt
  \belowdisplayshortskip=3pt plus 3pt 
  \smallskipamount=2pt plus 1pt minus 1pt
  \medskipamount=4pt plus 2pt minus 1pt
  \bigskipamount=9pt plus 3pt minus 3pt
  \normalbaselineskip=9pt
  \setbox\strutbox=\hbox{\vrule height7pt depth2pt width0pt}%
  \normalbaselines\rm}\catcode`\@=12
\newcount\noteno
\noteno=0
\def\up#1{\raise 1ex\hbox{\sevenrm#1}}
\def\note#1{\global\advance\noteno by1
\footnote{\parindent0.4cm\up{\number\noteno}\
}{\vtop{\eightpoint\baselineskip12pt\hsize15.5truecm\noindent
#1}}\parindent 0cm}
\vsize = 26truecm
\hsize = 16truecm
\voffset = -1truecm
\parindent=0cm
\baselineskip15pt
\overfullrule=0pt

\centerline{\bf The Coble hypersurfaces}
\smallskip
\smallskip \centerline{Arnaud {\pc BEAUVILLE}} 
\vskip1truecm
{\bf Introduction}\smallskip 
\ind The title refers to the following nice observations of
Coble. Let $A$ be a complex abelian variety, of dimension $g$, and
${\cal L}$ a line bundle on $A$ defining a principal polarization (that
is, ${\cal L}$ is ample and $\dim H^0(A,{\cal L})=1$). We will assume
throughout that $(A,{\cal L})$ is {\it indecomposable}, that is, cannot
be written as a product of principally polarized abelian varieties of
lower dimension.
\ind  Fix an integer
$\nu\ge 1$ and put
$V_\nu=H^0(A,{\cal L}^\nu)$. We consider the 
morphism\note{We use Grothendieck's notation: ${\bb  P}(V_\nu)$ is the space of hyperplanes of
$V_\nu$.}
$\varphi_\nu:A\rightarrow {\bb  P}(V_\nu)$  defined by the
global sections of ${\cal L}^\nu$. Recall that
$\varphi_\nu$ is an embedding for
$\nu\ge 3$, and that  $\varphi_2$
induces an embedding of the Kummer variety $A/\{\pm 1\}$ in 
${\bb  P}(V_2)$. Let  $A_\nu$ be the kernel of the multiplication by
 $\nu$ in
$A$; the  group
$A_\nu$ acts on $A$ and on 
${\bb  P}(V_\nu)$ in such a way that $\varphi_\nu$   is
$A_\nu$\tx equivariant. 
\smallskip 
{\bf  Proposition} (Coble)$.-$
 1) {\it  Let $g=2$. There exists a unique $A_3$\tx
invariant cubic hypersurface in ${\bb  P}(V_3)\ (\cong {\bb  P}^8)$ 
that is singular along  $\varphi_3(A)$. The polars of this
cubic span the space of quadrics in ${\bb  P}(V_3)$ containing}
$\varphi_3(A)$.

\ind 2) {\it  Let $g=3$. There exists a unique $A_2$\tx
invariant quartic hypersurface in ${\bb  P}(V_2)\ (\cong {\bb  P}^7)$ that is singular along 
$\varphi_2(A)$. The polars of this quartic span the space of
cubic hypersurfaces in ${\bb  P}(V_2)$ containing}
$\varphi_2(A)$.
  
\ind The proof of 2) appears in [C2], and that of 1) in [C1]
(actually the  cubic is not explicitly mentioned in that paper,
but it is easily deduced from the  equations for the
quadrics containing $\varphi_3(A)$. I am indebted to I.
Dolgachev for this reference). Both results are proved by
explicit computations. These hypersurfaces have a beautiful 
interpretation in terms of vector bundles on curves (see [N-R] for
the quartic and [O] for the cubic).
\ind An analogous statement appears in [O-P], this time for the
moduli space ${\cal SU}_C(2)$ of semi-stable rank 2 vector bundles
with trivial determinant on a non-hyperelliptic curve of genus 4
(this moduli space  is naturally embedded in ${\bb  P}(V_2)$). Oxbury
and Pauly prove that it is contained in a  unique $A_2$\tx invariant
quartic hypersurface, whose polars span the space of cubic
hypersurfaces containing 
${\cal SU}_C(2)$.
\ind The main observation of this note is that these facts follow from
a general (and elementary) result 
about representations of the Heisenberg group (Proposition 2.1
below). Let us just mention here a geometric
consequence of that result:\smallskip 
{\bf Proposition}$.-$
{\it Let $n=3$ or $4$; put $\nu=3$ if $n=3$, $\nu=2$
 if $n=4$. Let $(T_1,\ldots,T_N )$ be a coordinate system on ${\bb
P}(V_\nu)$. Let
$X$ be an
$A_\nu$\tx invariant subvariety of
${\bb  P}(V_\nu)$. 
Then the space of hypersurfaces of degree
$ n-1$ containing $X$ admits a basis $(\partial F_i/\partial
T_j)$, where $F_1,\ldots ,F_m$ are forms of degree $n$ on ${\bb
P}(V_\nu)$, such that the hypersurfaces $F_i=0$ are $A_\nu$\tx 
invariant {\rm (}and singular along}
$X$).
\vskip1truecm
{\bf 1. Heisenberg submodules of 
${\bf\sym^{\nn-1}V}$} 
\smallskip 
 \ind Let $n$ be an integer; we put $ \nu=n$ if $n$ is odd,
$\nu=n/2$ if
$n$ is even. We write for brevity $V$ instead of
$V_\nu$.  We will occasionally  pick a 
 coordinate system $(T_1,\ldots,T_N )$ on ${\bb P}(V_\nu)$, to make
some of our statements more concrete.
\ind The action of
$A_\nu$ on
${\bb  P}(V)$ lifts to an action on $V$ of  a central extension
$\widetilde{A}_\nu$ of
$A_\nu$ by ${\bb C}^*$.  For all $\gamma\in
\widetilde{A}_\nu$, the element $\gamma^{ n}$ belongs to the
center ${\bb C}^*$ of $\widetilde{A}_\nu$, and the map
$\gamma\mapsto\gamma^{ n}$ is a homomorphism of
$\widetilde{A}_\nu$ onto ${\bb C}^*$ (this is where we need to
take
$n=2\nu$ instead of $\nu$ when $\nu$ is even). We denote by
$H_n$ its kernel; it is a central extension
$$1\rightarrow \mug_{n}\longrightarrow H_n\longrightarrow
A_\nu\rightarrow 0\ $$of $A_\nu$ by the group $\mug_{n}$
of
${n}$\tx th roots of unity in ${\bb C}$.\smallskip  
{\bf Proposition 1}$.-$ 
 {\it Assume $n=3$ or $4$. Let $W$ be an irreducible sub-$\!H_n$\tx
module of 
$\sym^{ n-1}V$. There exists a  $H_n$\tx invariant
form
$F\in\sym^{ n}V$, unique up to a scalar, such that 
 $({\partial F/ \partial T_1},\ldots,{\partial F/ \partial T_N})$ form a
basis of $W$.}
\par
{\it Proof} :  Put $N=\dim V\ (=\nu^g)$. The group $H_n$ acts
irreducibly on $V$, and this is the unique irreducible
representation of $H_n$ on which the center 
$\mug_{n}$ acts by homotheties. It follows that 
  the representation of $H_n$ on $\sym^{ n-1}V$ is
isomorphic to the direct sum of  $k$  copies
of $V^*$, with
$$k=\dim \sym^{ n-1}V/\dim V^*={1\over N}{N+ n-2\choose  n-1}\
.$$  The space
$\Hom^{}_{H_n}(V^*,\sym^{ n-1}V)$ has dimension
$k$; it parametrizes the irreducible
sub-$\!H_n$\tx modules of $\sym^{ n-1}V$.
\ind   
 Consider the
$H_n$\tx equivariant injective map
$$h:\sym^{n} V\longrightarrow \Hom(V^*,\sym^{ n-1}V)$$
given by $h(F)(\partial)=\partial F$ (we identify $V^*$ with
the space of degree $-1$ derivations of $\sym V$).
It induces  an injection
$(\sym^{n} V)^{H_n}\mono \Hom_{H_n}(V^*,\sym^{ n-1}V)$ of the
 $H_n$\tx invariant subspaces. The assertion of the Proposition is
that this map is onto, or equivalently that $\dim (\sym^{n}
V)^{H_n}={1\over N}{N+ n-2\choose  n-1} $.
\ind    
The action of $H_n$ on $\sym^{n} V$ factors through the
abelian quotient $A_\nu$,  hence is the direct sum of $1$\tx
dimensional representations $V_\chi$ corresponding to characters
$\chi$ of $A_\nu$. We claim that all non-trivial characters of 
$A_\nu$ appear with the same multiplicity. To see this, consider
the group
$\Aut(H_n,\mug_n)$ of automorphisms of $H_n$ which induce the
identity on
$\mug_n$. Because of the unicity property of the representation
$\rho:H_n\rightarrow GL(V)$, for every $\varphi
\in\Aut(H_n,\mug_n)$ the representation
$\rho\rond\varphi $ is isomorphic to $\rho$, thus
$(\sym^n\rho)\rond\varphi $ is isomorphic to $\sym^n\rho$. 
This implies that the characters appearing in the decomposition of
$\sym^nV$ are exchanged by the action of $\Aut(H_n,\mug_n)$.
But the action of $\Aut(H_n,\mug_n)$ on
$A_\nu$ factors through 
 a surjective homomorphism $\Aut(H_n,\mug_n)\rightarrow {\rm
Sp}(A_\nu)$
  (see e.g. [B-L], ch. 6, lemma 6.6).  Since $\nu$ is prime, the
symplectic group ${\rm Sp}(A_\nu)$
 acts transitively on the set of
nontrivial characters of $A_\nu$, hence our claim.
\ind Thus we have
$$\sym^{n} V=
\bigl(\sdir_{\chi\not=1}^{}V_\chi\bigr)^m\oplus (\sym^{n}
V)^{H_n}$$ for some integer $m\ge 0$. Counting dimensions yields
$${N+ n-1\choose  n}=m\,(N^2-1)+\dim (\sym^{n} V)^{H_n} \ .$$ 
On the other hand a simple computation gives
$${N+ n-1\choose  n}=m\,(N^2-1)+{1\over N}{N+
n-2\choose  n-1} \ ,$$ with $m={1\over 6}(N+3)$ for $n=3$,
and $m={1\over 24}(N+2)(N+4)$ for $n=4$. Moreover we have 
$\dim (\sym^{n} V)^{H_n}\le {1\over N}{N+ n-2\choose 
n-1}<N^2-1$. Thus  $\dim (\sym^{n} V)^{H_n}$ and ${1\over
N}{N+ n-2\choose  n-1}$ are both equal to the rest of the division of
${N+ n-1\choose  n}$ by
$N^2-1$, hence they are equal.\cqfd
{\bf Corollary}$.-$ {\it Let $X$ be a
 subvariety of ${\bb  P}(V)$,  invariant under the action of
$A_\nu$; denote by ${\cal I}_X$  the ideal sheaf of $X$ in ${\bb 
P}(V)$. Let 
$(F_1,\ldots,F_m)$ be a basis of the space of
$H_n$\tx invariant forms
in $ \sym^{n}V$ which are singular along $X$. Then  the
partial derivatives 
 $({\partial F_i/ \partial T_j})$ form
a basis of  $H^0({\bb  P}(V),{\cal I}_X(n-1))$.
In particular, if $\dim H^0({\bb  P}(V),{\cal I}_X(n-1))=\nu^g$, there
exists a unique  $H_n$\tx invariant form
in $ \sym^{n}V$ which is singular along $X$.}
\ind Indeed $H^0({\bb  P}(V),{\cal I}_X(n-1))$ is a
sub-$\!H_n$\tx module of $H^0({\bb  P}(V),{\cal O}_{\bb
P}(n-1))=$ $\sym^{ n-1}V$, and therefore isomorphic to a direct sum  of
simple modules.\cqfd
\ind In the next section we will apply the Corollary to the abelian variety
$A$ embedded in ${\bb  P}(V_\nu)$. Another interesting case is
when $X$ is the moduli space of vector bundles of rank $2$ and trivial
determinant on a curve $C$ of genus $3$ with no vanishing
theta-constant.  Let $A$ be the Jacobian of $C$; then $X$ has a 
natural $A_2$\tx equivariant embedding in ${\bb P}(V_2)$, and
 Oxbury and Pauly prove the equality 
$\dim H^0({\bb P}(V_2),{\cal I}_X(3))=8$ [O-P]. Therefore there exists a
unique \hbox{$H_4$\tx invariant} quartic hypersurface singular along
$X$.
\medskip
 {\it Remark}$.-$  Unfortunately
the cases $n=3$ and $n=4$ seem to be the only ones for which
the Proposition holds. If for instance $n$ is prime $\ge 5$,
it is easy to check that  the equality $\dim (\sym^{n}
V)^{H_n}={1\over N}{N+n-2\choose n-1}$ {\it never} holds.
\smallskip 

\vskip1truecm
{\bf 2. Application: equations for abelian varieties}\smallskip 
\ind (2.1) Let us apply the Corollary to $X=\varphi _\nu(A)$
embedded in ${\bb  P}(V_\nu)$.  If
$n=4$ we will assume that
$(A,{\cal L})$ {\it has no vanishing theta-constant} (that is, no symmetric 
theta divisor singular at $0$ -- if $g=3$ this simply means that $(A,{\cal
L})$ is the Jacobian of a non-hyperelliptic curve). This implies that the
Kummer variety
$\varphi_2(A)\i{\bb  P}(V_2)$ is projectively normal, while
$\varphi_3(A)$ is always projectively normal in
${\bb  P}(V_3)$ [Ko]. Thus the natural map
$H^0({\bb  P}(V_\nu),{\cal O}_{\bb  P}(n-1)) \rightarrow H^0(X,
{\cal O}_X(n-1))$ is surjective, and this allows us to compute the
dimension of its kernel. We find that {\it
 the space of $H_n$\tx invariant forms
in $ \sym^{n}V$  singular along $X$ has dimension $m_n(g)$ given
by}
$$m_3(g)={1\over 2}(3^g-2^{g+1}+1)
\qquad  m_4(g)={1\over
6}(2^g(2^g+3)-3^{g+1}-1)\
;$$for any basis $(F_1,\ldots,F_{m_n(g)})$ of this space, the
derivatives $ (\partial F_i/ \partial T_j)$ form a
basis of the space
of forms of degree $n-1$ vanishing along $X$. 
\smallskip 
\ind (2.2) Let us consider in particular the case $g=n-1$ considered
by Coble. Since
$m_3(2)=m_4(3)=1$ we recover Coble's result: there is a unique
$H_n$\tx invariant  hypersurface  of degree $n$
singular along $\varphi_\nu(A)$. In fact  we have a slightly
better result:\smallskip {\bf Proposition 2}$.-$ {\it Assume $g=n-1$.
The Coble hypersurface in ${\bb  P}(V_\nu)$ is the unique
hypersurface of degree $n$ singular along $\varphi_\nu(A)$.}
\par
{\it Proof} : The case of the Coble quartic is explained in [L], and the
proof works equally well for the cubic. Let us recall briefly the
argument. Let $F=0$ be the Coble hyper\-surface. The derivatives
${\partial F/ \partial T_1},\ldots,{\partial F/ \partial T_N}$ span the
space $I_{n-1}$ of forms of degree $n-1$ vanishing along
$\varphi_\nu(A)$;  the action of $H_n$ on $I_{n-1}$ is  irreducible.
\ind Let $W$ be the space 
of forms of degree $n$ which are singular along
$\varphi_\nu(A)$; it is a sub-$\!H_n$\tx module of $\sym^nV$,
hence a sum of one-dimensional  representations $W_\chi$. Let 
$G\not=0$ in
$W_\chi$. The derivatives
${\partial G/ \partial T_1},\ldots,{\partial G/ \partial T_N}$ vanish
on $\varphi_\nu(A)$, hence span a subspace of $I_{n-1}$; since
this subspace is stable under $H_n$, it is equal to   $I_{n-1}$. By [D,
\S 1] this implies that there exists an automorphism $T$ of $V_\nu$ such that $G=F\rond T$. 
\ind Now the singular locus of the Coble hypersurface is exactly
$\varphi_\nu(A)$ (see (2.3) below); thus $T$ must preserve
$\varphi_\nu(A)$. In the group of automorphisms of
$V_\nu$ preserving $\varphi_\nu(A)$, the Heisenberg group
$H_n$ is normal --  because the group of translations of $A$ is
normal inside the group of all automorphisms. Thus $T$ normalizes
$H_n$; this implies that the form $G=F\rond T$ is $H_n$\tx
invariant, and therefore proportional to $F$ by Coble's result.\cqfd
\ind (2.3) For $g=2$, Coble
states in [C1] that
$\varphi _3(A)$  is the set-theoretical intersection of the quadrics 
that contain it -- in other words, $\varphi _3(A)$ is the singular locus
of the Coble cubic; this is proved even scheme-theoretically in [B].
When
$g=3$ and $(A,{\cal L})$ has no vanishing theta-constant,
Narasimhan and Ramanan have proved that the Kummer variety
$\varphi _2(A)$ is set-theoretically the singular locus of the Coble
quartic [N-R]; this holds also scheme-theoretically by  [L]. It is tempting
to conjecture that  both statements hold in higher dimension as well,
namely that the  abelian variety $\varphi _3(A)$ is a 
scheme-theoretical intersection of  quadrics and that the Kummer
variety $\varphi _2(A)$ is a 
scheme-theoretical intersection of  cubics. Note, however, that these
quadrics or cubics cannot generate the full ideal of $\varphi_\nu(A)$:
\smallskip 
{\bf Proposition 3}$.-$ {\it The graded ideal $I$ of 
$\varphi_\nu(A)$ in
${\bb  P}(V_\nu)$ is} not {\it generated by its
elements of degree $\le n-1$}.\ind 
(Recall that $I$ is generated by its
elements of degree $\le n$, see [B-L], ch. 7 and [K].)
\ind Note that  the Proposition
 is immediate in the case $g=n-1$ considered by Coble, 
because  then
$\dim(V\otimes I_{n-1})<\dim I_n$. However  this inequality does
not hold any more in higher genus.\smallskip 
\par {\it Proof}:
 We will prove the inequality    $\dim (V\otimes
I_{n-1})^{H_n}<\dim (I_{n})^{H_n}$, which implies that the
multiplication map $V\otimes I_{n-1}\rightarrow I_n$ cannot be
surjective. Let us treat first the case $n=3$. From the exact sequence
$0\rightarrow I_3\rightarrow\sym^3V\rightarrow H^0(A,{\cal
L}^{9})\rightarrow 0$ (2.1) we get $$\dim
I_3={N+2\choose 3}-N^2={N-3\over 6}(N^2-1)+{N-1\over 2} \ ;$$
as in Proposition 1  we conclude that $\dim (I_3)^{H_3}=(N-1)/2$. 
\ind Let $K\i\sym^3V$ be the space of $H_3$\tx invariant
cubic forms singular along $\varphi_3(A)$; by the
proposition the natural map $V^*\otimes K\rightarrow I_2$ is
an isomorphism. The action of
$H_3$ on $K$ is trivial, and the $H_3$\tx module $V\otimes
V^*$ is the direct sum  of a one-dimensional factor for each
character of $A_3$; thus $$\dim (V\otimes I_2)^{H_3}=\dim
K={1\over 2}(3^g-2^{g+1}+1)<{1\over 2}(N-1)$$(2.1),
hence the result.
\ind For $n=4$ the same method gives $\dim
(I_4)^{H_2}={1\over 6}(N-1)(N-2)$, which is larger than
 $\dim (V\otimes
I_3)^{H_2}={1\over
6}(N(N+3)-3^{g+1}-1)$.\cqfd
\vskip2truecm
\centerline{ REFERENCES} \vglue15pt{\baselineskip12.5pt
\def\num#1{\smallskip\item{\hbox to\parindent{\enskip [#1]\hfill}}}
\parindent=1.3cm 
\num{B} W. {\pc BARTH}: {\sl Quadratic equations for
level-{\rm 3} abelian surfaces}. Abelian varieties, 1--18; de
Gruyter, Berlin -- New York (1995).
\num{B-L} C. {\pc BIRKENHAKE}, H. {\pc LANGE}: {\sl Complex
abelian varieties}. Grund.
   Math. Wiss. {\bf 302}, Springer-Verlag,
   Berlin, 1992.
 \num{C1} A. {\pc COBLE}: {\sl Point Sets and Allied Cremona
Groups} III. Trans.  Amer. Math. Soc. {\bf  18} (1917),  331--372. 
\num{C2} A. {\pc COBLE}: {\sl Algebraic geometry and theta functions}.
Amer. Math. Soc. Colloquium Publi. 10 (1929). Amer. Math. Soc.,
Providence (1982).
\num{D} R. {\pc DONAGI}: {\sl Generic Torelli for projective
hypersurfaces}. Compositio Math. {\bf 50} (1983),  325--353.
\num{K} A. {\pc KHALED}: {\sl \'Equations des vari\'et\'es de
Kummer}.  Math. Ann. {\bf 295}
(1993),  685--701.
\num{Ko} S. {\pc KOIZUMI}: {\sl Theta relations and projective
normality of Abelian varieties}. Amer. J. Math. {\bf 98} (1976), 
865--889.
\num{L} Y. {\pc LASZLO}: {\sl Local structure of the moduli space
of vector bundles over curves}. Comment. Math. Helv. {\bf 71}
(1996),  373--401. 
\num{N-R} M.S. {\pc NARASIMHAN}, S. {\pc RAMANAN}: {\sl
$2\theta$\tx linear  systems  on Abelian varieties.} Vector bundles
on algebraic varieties, 415--427; Oxford University Press (1987).
\num{O-P} W. {\pc OXBURY}, C. {\pc PAULY}: {\sl Heisenberg invariant 
quartics and ${\cal SU}_C(2)$ for a curve of
   genus four}. Math. Proc. Cambridge Philos. Soc. {\bf  125} (1999),
295--319. 
\vskip1cm
\def\pc#1{\eightrm#1\sixrm}
\hfill\vtop{\eightrm\hbox to 5cm{\hfill Arnaud {\pc BEAUVILLE}\hfill}
\vskip2pt
 \hbox to 5cm{\hfill Institut Universitaire de France\hfill}\vskip-4pt
\hbox to 5cm{\hfill \&\hfill}\vskip-4pt
 \hbox to 5cm{\hfill Laboratoire J.-A. Dieudonn\'e\hfill}
 \hbox to 5cm{\sixrm\hfill UMR 6621 du CNRS\hfill}
\hbox to 5cm{\hfill {\pc UNIVERSIT\'E DE}  {\pc NICE}\hfill}
\hbox to 5cm{\hfill  Parc Valrose\hfill}
\hbox to 5cm{\hfill F-06108 {\pc NICE} Cedex 2\hfill}}}
\end